\documentclass{amsart}
\usepackage{amscd,graphicx}

\DeclareMathSymbol{\twoheadrightarrow}  {\mathrel}{AMSa}{"10} \newcommand{\onto}{\twoheadrightarrow}

\newcommand{\Gss}{G^{\rm ss}}

\newcommand{\inv}{^{-1}}

 \newcommand{\diag}{\operatorname{diag}}

\newcommand{\mld}{\operatorname{mld}}
\newcommand{\im}{\operatorname{im}}

\newcommand{\GL}{\operatorname{GL}}

\newcommand{\chr}{\operatorname{char}}
\newcommand{\supp}{\operatorname{supp}}
\newcommand{\Hom}{\operatorname{Hom}}
\newcommand{\Aut}{\operatorname{Aut}}

\newcommand{\hcD}{{\widehat{\mathcal D}}} \newcommand{\wcD}{{\widetilde{\mathcal D}}}

\newcommand{\bN}{{\mathbb N}}
\newcommand{\bG}{{\mathbb G}}
\newcommand{\bQ}{{\mathbb Q}}

\newcommand{\bZ}{{\mathbb Z}}

\newcommand{\cD}{{\mathcal D}}

\newcommand{\cV}{{\mathcal V}}
\newcommand{\cX}{{\mathcal X}}
\newcommand{\cS}{{\mathcal S}}

\newcommand{\wX}{{\widetilde X}}
\newcommand{\wD}{{\widetilde D}}
\newcommand{\wN}{{\widetilde N}}
\newcommand{\wH}{{\widetilde H}}
\newcommand{\tv}{{\tilde v}}
\newcommand{\td}{{\tilde d}}
\newcommand{\tf}{{\tilde f}} \newcommand{\wcV}{{\widetilde{\mathcal V}}}

\newcommand{\nq}{N_{\bQ}}

\theoremstyle{plain}
\newtheorem{theorem}{Theorem}[section]

 \newtheorem{lemma}[theorem]{Lemma}
\newtheorem{corollary}[theorem]{Corollary}

\newtheorem{conjecture}[theorem]{Conjecture}

\theoremstyle{definition} \newtheorem{definition}[theorem]{Definition}

\newtheorem{notation}[theorem]{Notation}         

\newtheorem{example}[theorem]{Example}

\newtheorem{remark}[theorem]{Remark}

\newtheorem*{acknowledgements-no}{Acknowledgments} 

\theoremstyle{remark}

\author{Valery Alexeev and Michel Brion}
\title{Boundedness of spherical Fano varieties}
\date{January 17, 2003}

\begin{document}
\bibliographystyle{amsalpha}
\maketitle

Classically, G. Fano proved that the family of (smooth,
anticanonically embedded) Fano 3-dimensional varieties is bounded, and
moreover provided their classification, later completed by V.A.
Iskovskikh, S.  Mukai and S. Mori. For singular Fano varieties with
log terminal singularities, there are two basic boundedness
conjectures: Index Boundedness and the much stronger $\epsilon$-lt
Boundedness.

The $\epsilon$-lt Boundedness was known only in two cases: in
dimension 2 \cite{Alexeev_Boundedness} and for toric varieties
\cite{Borisovs}. In this paper we prove it for a significantly less
``elementary'' class, that of spherical varieties. In addition to an
argument adapted from the toric case, the proof contains quite a few
new twists.

In Section~\ref{sec:Fano invariant of a spherical subgroup}, we
introduce a new invariant of a spherical subgroup $H$ in a reductive
group $G$ which measures how nice an equivariant Fano compactification
of $G/H$ there exists.

\section{Boundedness conjectures}

Let $X$ be a log terminal Fano variety defined over an algebraically
closed field $k$,~i.e.
\begin{enumerate}
\item $X$ is normal,
\item the canonical class $K_X$ is $\bQ$-Cartier,
\item for a log resolution of singularities $f:Y\to X$, in the formula
  $$
  K_Y = f^*K_X + \sum a_iE_i $$
  the log discrepancies $b_i = 1+
  a_i$ of exceptional divisors $E_i$ are all positive, \item $-K_X$ is
  ample. \end{enumerate}

\begin{definition}
  The \emph{index} of $X$ is the minimal positive integer $I$ such
  that $I\cdot K_X$ is Cartier.
\end{definition}

\begin{definition}
  The \emph{minimal log discrepancy} of $X$ is $\mld(X) = \min b_i$.
  For $\epsilon>0$, we will say that $X$ is $\epsilon$-lt if
  $\mld(X)>\epsilon$.
\end{definition}

\begin{definition}
  One says that a certain set of varieties (or schemes) is
  \emph{bounded} if they appear as geometric fibers of a family
  $\cX\to \cS$, with $\cX$ and $\cS$ both schemes of finite type over
  the base field.
\end{definition}

The following conjecture appeared independently in
\cite{Alexeev_Boundedness} and the work of Alexander and Lev Borisovs
\cite{Borisovs}. \begin{conjecture}[$\epsilon$-lt Boundedness or BAB
  Conjecture] The family of $n$-dimensional $\epsilon$-lt Fano
  varieties is bounded.
\end{conjecture}
 
The only cases where this conjecture has been proven are:
\begin{enumerate} \item $n=2$, arbitrary $\chr k$
  \cite{Alexeev_Boundedness}, \item for toric Fano varieties,
  arbitrary $\chr k$ \cite{Borisovs}, \item $n=3$ and $\epsilon=1$
  (i.e. with terminal singularities), $\chr k=0$
  \cite{Kawamata_Boundedness}.
\end{enumerate}

A weaker conjecture (which follows from the previous one by taking
$\epsilon=1/I$) is due to V. Batyrev: \begin{conjecture}[Index
  Boundedness] The family of $n$-dimensional log terminal Fano
  varieties of index $I$ is bounded.
\end{conjecture}
This conjecture is known in the following cases not covered by the
above ($\chr k=0$ everywhere): \begin{enumerate} \item for smooth
  Fanos \cite{KollarMiyaokaMori92b}, \item $n=3$
  \cite{Borisov_3foldBoundedness}. \end{enumerate} Recently, J.
McKernan \cite{McKernan_Boundedness} announced a proof of Index
Boundedness Conjecture for any $n$.

For the rest of the paper, we will be working over an algebraically
closed field of characteristic $0$.

\section{Recall of spherical varieties}
\label{Recall of spherical varieties}

Let $T=\bG_m^n$ be a multiplicative torus. Toric varieties are of
course normal varieties with an open $T$-orbit which are (partial)
compactifications of $T$.  They have a familiar combinatorial
description. One starts with two free abelian groups of rank $n$, the
group of characters $M$ of $T$ and the dual group $N$. Every toric
variety $X$ corresponds uniquely to a fan $\Sigma$ in $\nq$, and $X$
is complete iff $\supp\Sigma = \nq$.

$T$-invariant Cartier divisors correspond to piecewise linear
functions on the fan. Let $v\in\Sigma(1)$ be the integral generators
of 1-dimensional cones, they correspond to $T$-invariant divisors
$D_v$.  Then $-K_X = \sum D_v$. In particular, $K_X$ is
$\bQ$-Gorenstein iff for any cone $\sigma$ its integral generators
$v\in \sigma(1)$ lie on a common hyperplane.

Spherical varieties generalize this picture to the non-commutative
case. Let $G$ be a connected reductive group with a Borel subgroup
$B$. A variety $X$ with $G$-action is called \emph{spherical} if it is
normal and has an open $B$-orbit.  If $H$ is the stabilizer of a point
$x$ in the open orbit then $X$ is a partial compactification of $G/H$.
Luna and Vust proved the structure theory of spherical varieties which
was later translated into the language of colored fans.  The basic
references for this theory are \cite{Knop_LV} and \cite{Brion_Notes}.
Note that the class of spherical varieties is much richer than that of
toric varieties, and they may have fairly complicated singularities
(for example $\bQ$-factorial spherical singularities need not be
quotient).

\begin{remark}\label{caveat}
  A major caveat is that the whole theory is relative to $G/H$. Once
  the homogeneous variety $G/H$ is fixed, spherical compactifications
  correspond uniquely to colored fans.  The question of which
  subgroups $H\subset G$ are spherical is open, with complete
  classification available only for groups of type A \cite{Luna}. In
  any case, for any noncommutative group $G$ there exist infinitely
  many non-isomorphic spherical subgroups, and infinitely many
  non-isomorphic colored lattices corresponding to them. \end{remark}

Let us summarize all the combinatorial facts that we need and fix the
notation. Let $X$ be a spherical embedding of $G/H$. Then we have:
\begin{enumerate}
\item a free abelian group $\Lambda$ of rank $r\le n$ (analog of $M$
  for toric varieties), it is defined as the multiplicative group of
  $B$-eigenvectors of $k(G/H)$ up to scalars,
  $\Lambda=k(G/H)^{(B)}/k^*$;
\item its dual $\Lambda^*= \Hom(\Lambda,\bZ)$, which we will denote by
  $N$, by analogy with the toric case;
\item the vector space $\nq= \Lambda^*_{\bQ}$.
\end{enumerate}
The new elements are:
\begin{enumerate}\setcounter{enumi}{3}
\item a \emph{valuation cone} $\cV$; this is a co-simplicial cone
  spanning $\nq$;
\item a finite set $\cD$ of \emph{colors} (these are the irreducible
  $B$-stable divisors in $G/H$) and a map $\rho:\cD\to N$.
\end{enumerate}

A simple embedding, i.e. one with a unique closed $G$-orbit $Y$,
defines two finite sets $\cV_Y\subseteq \cV$ and $\cD_Y\subseteq \cD$.
The corresponding colored cone is the pair $(\sigma_Y, \cD_Y)$ where
the cone $\sigma_Y$ is generated by $\cV_Y$ and $\rho(\cD_Y)$. It is
known that $\sigma_Y$ is strictly convex and that
$\sigma_Y^0\cap\cV\ne\emptyset$.  Variety $X$ is uniquely described by
a colored fan $\Sigma$, i.e.  a compatible collection of colored
cones. $X$ is complete iff the cones cover $\cV$; in this case the
collection $\{\sigma\cap\cV\}$ is an ordinary fan with support $\cV$.
In particular, the embeddings without colors correspond to partial
subdivisions of $\cV$. Such embeddings are called \emph{toroidal}.

An element of $\cV_X = \cup \cV_Y$, resp. $\cD_X= \cup \cD_Y$, where
the sum goes over the finite set of $G$-orbits $Y$, defines a
$G$-invariant, resp. $B$-invariant but not $G$-invariant, divisor on
$X$. Every $B$-invariant divisor can be written as
$$
(\sum_{v\in \cV_X} n_v D_v + \sum_{D\in \cD_X} n_D D) + \sum_{D\in
  \cD\setminus\cD_X}n_D D.
$$
This divisor is Cartier iff the coefficients in the first half of
this expression are values of a piecewise linear function
$(l_{\sigma}\in \Lambda, \sigma\in \Sigma)$ at the points $v$,
$\rho(D)$.  The divisors $D\in \cD\setminus\cD_X$ are always Cartier.
By \cite{Brion_CurvesDivs}, the anticanonical divisor on a spherical
variety can be written as \begin{eqnarray}\label{eq:canclass} -K_X =
  \sum_{v\in \cV_X} D_v + \sum_{D\in \cD} a_D D \end{eqnarray} for
unique, explicitly computable positive integers $a_D$.  Variety $X$ is
$\bQ$-Gorenstein iff for every colored cone $(\sigma_Y, \cD_Y)$ the
points $v\in\sigma_Y(1)$ and $\rho(D)/a_D$, $D\in \cD_Y$, lie on a
common hyperplane.

\section{Boundedness of toric Fanos}

Let us first recall the proof of boundedness in the toric case. The
anticanonical divisor on a $\bQ$-Gorenstein toric variety corresponds
to a piecewise linear function $l$ taking value 1 at every integral
generator $v\in \Sigma(1)$.  This divisor is ample iff the function is
strictly convex, i.e. if $v$'s are vertices of a convex lattice
polytope, call it $Q$. The log discrepancies at $T$-invariant
exceptional divisors are precisely the values of $l$ at integral point
$u\in N$. Since every toric variety has a $T$-equivariant resolution,
$X$ is $\epsilon$-lt if and only if the following condition is
satisfied (where $Q^0$ denotes the interior): \begin{eqnarray}
  \label{eq:elt}
\epsilon Q^0\cap N = \{0\}  
\end{eqnarray}

Hence, the boundedness follows from the following theorem,
\cite{Hensley}:
\begin{theorem}[Hensley]\label{thm:Hensley}
  Up to $\Aut(N)=\GL(n,\bZ)$, there are only finitely many lattice
  polytopes in $\nq$ satisfying the condition (\ref{eq:elt}).
\end{theorem} Borisovs gave a different proof of Theorem
\ref{thm:Hensley} in \cite{Borisovs}.

\section{Boundedness of spherical Fanos}

\begin{notation}
  Let $X$ be a spherical $G$-variety with open orbit $G/H$.  Denote by
  $C=Z(G)^0$ the connected component of the center, and by
  $\Gss=[G,G]$ the derived subgroup, the semisimple part of $G$.  Then
  $G$ is the quotient of $\Gss \times C$ by a finite central subgroup.
  We will also denote by $\wH=N_G(H)$ the normalizer of $H$ in $G$.
  Then $\wH$ contains $C$, and $\wH/H=\Aut^G(G/H)$, the automorphism
  group of the $G$-variety $G/H$. In fact, any element of $\wH/H$
  extends to an automorphism of $X$, i.e., $\Aut^G(X)=\wH/H$. Further,
  this group is diagonalizable.
\end{notation}

Here is our main theorem:

\begin{theorem}\label{thm:main}
  For any $\epsilon>0$ and $n\in \bN$, the set of $\epsilon$-lt
  $n$-dimensional spherical Fano varieties is finite. \end{theorem} As
a first reduction, we can assume that the action is faithful (just
divide by the kernel).  However, it will be convenient to work with
\emph{almost faithful} actions, in the following sense:
\begin{definition} A $G$-action on $X$ is called \emph{almost
    faithful} if its kernel is finite, and $C$ acts faithfully.
\end{definition}
As a second reduction, we can assume that $G = \Gss\times C$ with
simply connected $\Gss$ (go to a finite cover); then the Picard group
of $G$ is trivial. We can also assume that the action is almost
faithful (divide by a finite subgroup of $C$, if necessary).

\begin{proof}[Proof of Theorem~\ref{thm:main}]
  By Theorem~\ref{thm:refined} below, it suffices to bound the
  dimension of $G$ in terms of $n=\dim(X)$, where $X$ is
  $G$-spherical. This will imply that only finitely many isomorphism
  classes of connected reductive groups $G$ occur, and our theorem
  will follow.
  
  The dimension of $C$ is bounded by $n$, of course.  To bound the
  dimension of $\Gss$, consider a toroidal resolution $\wX$ of $X$.
  All closed orbits in $\wX$ are isomorphic to the same flag variety
  $G/P=\Gss/(P\cap \Gss)$, on which $\Gss$ acts with a finite kernel.
  And, by a result of Akhiezer \cite{Akhiezer}, the dimension of $\Aut
  G/P$ is bounded by a function of $\dim G/P \le n$ only. \end{proof}

\begin{definition}
  A $G$-action on $X$ is called \emph{smart} if
  \begin{enumerate}
  \item the action is almost faithful, and
  \item the natural homomorphism $C\to \Aut^G(X)^0$ is an isomorphism.
  \end{enumerate}
\end{definition}

\begin{lemma}\label{lem:reductions}
  Any almost faithful $G$-action can be extended to a smart action of
  a bigger connected reductive group.
\end{lemma}
\begin{proof}
  Consider the map
  $$
  \varphi:G=\Gss\times C \to \tilde G = \Gss \times (\wH/H)^0,
  ~(g,c)\mapsto (g,cH). $$
  Then $\varphi$ is a group homomorphism,
  injective since $C\cap H$ is trivial (as $C$ acts faithfully).
  Further, $\tilde G$ acts on $G/H$ by $(g,\gamma)\cdot xH = gx\gamma
  H$, this action extends to $X$, and $\varphi$ is equivariant. We
  check that the $\tilde G$-action is smart.
  
  Clearly, the derived subgroup $\Gss$ acts with finite kernel. Thus,
  the connected kernel of the $\tilde G$--action is contained in the
  connected center $(\wH/H)^0$.  But the latter acts faithfully, so
  that (1) holds. Finally, since $G$ embeds into $\tilde G$, we have
  $\Aut^{\tilde G}(X)^0\subseteq \Aut^G(X)^0 = (\wH/H)^0$ which
  implies (2). \end{proof}

\begin{definition}
  Two $G$-actions are called \emph{equivalent} if they differ by an
  automorphism of $G$.
\end{definition}
If $G=\Gss\times C$ then $\Aut G$ contains as a subgroup of finite
index the product of the group of inner automorphisms of $\Gss$ and of
$\Aut C \simeq \GL(\dim C, \bZ)$.

\begin{theorem}\label{thm:refined}
  For a fixed connected reductive group $G=\Gss\times C$ and
  $\epsilon>0$, there exist only finitely many spherical $G$-varieties
  with smart action which are $\epsilon$-lt Fano varieties, up to
  choosing an equivalent action.
\end{theorem}

To prepare for the proof, let us first translate in terms of colored
fans the conditions that $X$ is Fano and $\epsilon$-lt.  Let $X$ be
our spherical $G$-variety.  Therefore, we have a lattice $N$,
valuation cone $\cV$ and colors $\rho:\cD\to N$.  This time, the
anticanonical divisor corresponds to a piecewise linear function
$l=(l_{\sigma}, \sigma\in \Sigma)$ which takes values~1 at the points
$v\in \cV_X$ and $\rho(D)/a_D$.  Ampleness translates into two
conditions (see \cite[Ch.5.2]{Brion_Notes}):
  \begin{enumerate}
  \item $l$ is strictly convex; in other words $l_{\sigma} <
    l_{\sigma'} $ on $\sigma'\setminus \sigma'\cap\sigma$, and
  \item $l_{\sigma}(\rho(D)/a_D) = 1$ for every $D\in \cD_{\sigma}$,
    and $<1$ for every $D\in \cD\setminus \cD_{\sigma}$.
  \end{enumerate}
  This means that the linear inequalities $l_{\sigma} \le 1$ define a
  possibly unbounded polyhedral body, call it $P$, such that
  \begin{enumerate}
  \item the points $v\in \cV_X$ and $\rho(D)/a_D$, $D\in \cD_X$, are
    vertices of $P\cap |\Sigma|$, and
  \item the points $\rho(D)/a_D$, $D\in \cD\setminus\cD_X$, are in the
    interior $P^0$.
  \end{enumerate}
  
  Every spherical variety has a toroidal $G$-equivariant resolution,
  and the log discrepancies are precisely the values of $l$ at
  integral points $u\in N\cap\cV$.  Therefore, the $\epsilon$-lt
  condition is equivalent to
  \begin{eqnarray}
    \label{eq:cond_onP}
  \epsilon P^0 \cap N \cap \cV = \{0\}
  \end{eqnarray}
  Now consider the convex hull $Q$ of all points $v\in \cV_X$ and
  $\rho(D)/a_D$.  It is known that the cone $\cV$ and the set
  $\{\rho(D)/a_D\}$ do not lie in a common half-space of $\nq$, see
  \cite[Rem.3.4(3)]{Brion_Notes}. This implies that $Q$ is a polytope
  of maximal dimension whose interior contains the origin. Clearly,
  $Q$ also satisfies the condition (\ref{eq:cond_onP}) above.

  The role of Hensley's theorem will be played by the following
  combinatorial statement: \begin{lemma}
  \label{lem:comb_finiteness}
  Suppose we have fixed the following data:
  \begin{enumerate}
  \item a surjective homomorphism of lattices $\pi:N\to \wN_1$,
  \item a strictly convex rational polyhedral cone $\wcV$ generating
    $\wN_1\otimes\bQ$ and a cone $\cV=\pi\inv(\wcV)$,
  \item a positive integer $A$,
  \item finitely many points $\td_i=\pi(d_i)\in \wN_1/A$.
  \end{enumerate}
  Then for every $\epsilon>0$ there exist only finitely many
  maximal-dimensional polytopes $Q$ in $\nq$ satisfying the following
  condition:
    \begin{enumerate}
    \item $Q$ is the convex hull of some $d_i\in N/A$ with $\pi(d_i) =
      \td_i$ and finitely many elements of $\cV\cap N$,
    \item $\epsilon Q^0 \cap N \cap \cV = \{0\}$;
    \end{enumerate}
    up to the action of the group $\Aut(N,\pi)\subseteq \Aut(N)$
    leaving $\ker\pi$ invariant and inducing the identity on $\wN_1$.
  \end{lemma}

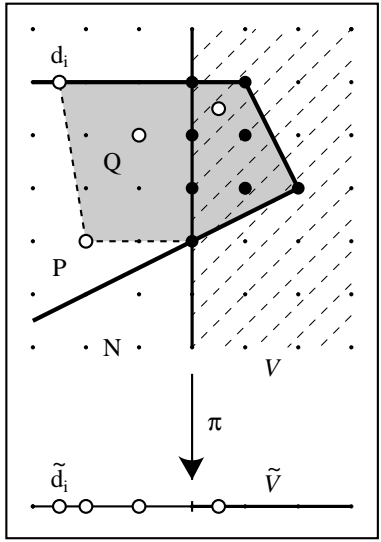
\begin{figure}[htbp]
  \includegraphics[scale=.8]{c-color.eps}
\end{figure}
\begin{proof}
  Of course, by replacing lattices $N,\wN_1$ by $N/A,\wN_1/A$ we can
  assume from the start that $Q$ is a lattice polytope. Polyhedral
  cone $\cV$ is defined by several linear inequalities $f_i(x)\ge 0$
  in $\nq$, where $f_i=\tf_i\circ\pi$ for some linear integral-valued
  functions on $\wN_1$.  Let $-F= \min \tf_j(\td_i)$.  For every point
  $x$ of the polytope $Q$ one has $f_j(x)\ge -F$. Therefore, if
  $\epsilon<1/F$ then
  $$
  \epsilon Q^0 \cap N \cap \cV = \epsilon Q^0 \cap N =\{0\}.
  $$
  Therefore, by Hensley's theorem there exist only finitely many
  polytopes satisfying our condition, up to $\Aut(N)$. Now let us fix
  one of those polytopes $Q$ and ask: for which $g\in\Aut(N)$ does the
  polytope $gQ$ has the same shape -- it is the convex hull of some
  points $d_i$ with $\pi(d_i)=\td_i$ and some elements in $\cV\cap N$?
  
  The vertices of $Q$ split into two sets: vertices $v_j$ that are in
  $\cV$, and vertices lying outside of $\cV$. The latter are
  necessarily some of $d_i$'s, say $d_1\dots d_s$. Up to a permutation
  of vertices (finitely many choices), we can assume that $\pi(g
  d_i)=\pi(d_i)$, $i=1\dots s$, and that $g v_j\in\cV$. Since $0\in
  Q$, there exist some positive integers $n_i,m_j$ such that
  $$
  \sum m_j v_j + \sum_{i=1}^s n_i d_i =0.
  $$
  Therefore, for some positive combination of $v_j$, one has
  $$
  \tv = \sum m_j \pi(v_j) = \sum m_j \pi(gv_j).
  $$
  Consequently, each of $gv_j$ belongs to $\cV + (\pi\inv(\tv)
  -\cV)$, and since $\wcV$ is strictly convex, this is a finite union
  of cosets of the lattice $\ker\pi$. Hence, up to finitely many
  choices again, we can assume that $gv_j\in v_j + \ker \pi$. Together
  with $gd_i\in d_i +\ker\pi$, $i=1\dots s$, this implies $g\in
  \Aut(N,\pi)$.
\end{proof}

We will also use the following basic boundedness result from
\cite{AlexeevBrion_Moduli} which follows from Knop's theorem about
local rigidity of wonderful compactifications:
\begin{theorem}\label{thm:bound_normalized} For any connected
  reductive group $G$, up to conjugation, there exist only finitely
  many spherical subgroups $H$ (i.e. $G/H$ is spherical) which
  coincide with their normalizer. \end{theorem}

We are now ready for our theorem.

\begin{proof}[Proof of Thm. \ref{thm:refined}]
  Let $X$ be an embedding of $G/H$.  By
  Theorem~\ref{thm:bound_normalized}, we can fix the normalizer $\wH$
  from now on. For the connected components of identity we also have
  $H_0\subseteq \wH_0$. We now have four lattices, $N^*=
  k(G/H)^{(B)}/k^*$ etc., of rank bounded by the rank of $G$, and a
  diagram
  $$
\begin{CD}
  N_0 @>>> N \\
  @V\pi_0VV     @VV\pi V \\
  \wN_0 @>>> \im N @>>> \wN.
\end{CD}
$$
In this diagram, horizontal arrows are injective with finite
cokernels, and vertical arrows are surjective ($\pi_0$ is surjective
because $H_0,\wH_0$ are connected). The map $N\to \wN$ need not be
surjective. However, there are only finitely many possibilities for
intermediate lattice $\wN_1=\im N$, so we can fix one of them.

The sets of colors $\cD,\wcD$ for groups $H,\wH$ are nearly the same.
More precisely, let $B$ be a Borel subgroup such that $BH$ is open in
$G$, then the colors of $G/H$ are the irreducible components of $(G
\setminus BH)/H$. But $BH=B\wH$, so that the set $\hcD$ of irreducible
components of $G\setminus BH = G \setminus B\wH$ is the same. The sets
of colors are obtained from $\hcD$ by dividing by $H$, resp. $\wH$,
which identifies some of the divisors into one color. We have a
commutative diagram
$$
\begin{matrix}
  \hcD & \onto & \cD & \overset{\rho}{\longrightarrow} N \\
  || & & \downarrow & {\hskip 25pt}\downarrow \pi \\
  \hcD & \onto & \wcD & \overset{\tilde{\rho}}{\longrightarrow} \wN_1.
\end{matrix}
$$
We will denote by $d_i$ (resp. $\td_i$) the points $\rho(D)/a_D$,
$D\in \cD$, (resp.  $\tilde\rho(\wD)/a_D$, $\wD\in \wcD$) and by $A$
the $\operatorname{LCM}(a_D)$.  By Lemma~\ref{lem:comb_finiteness},
there are only finitely many possibilities for the polytope $Q$ in
$\nq$, and hence also for the set $\rho(\cD)$, up to the action of
group $\Aut(N,\pi)$.

Consider the multiplicative group $k(G)^{(B\times H)}$ of those
rational functions on $G$ that are eigenvectors of $B$ (for left
multiplication) and of $H$ (for right multiplication).  Any such
function may have zeroes and poles along elements of $\hcD$ only;
conversely, any divisor in $\hcD$ has a global equation (since $G$ has
trivial Picard group) which is a $B\times H$-eigenvector. Moreover,
the functions in $k(G)^{(B\times H)}$ having no zero or poles are the
regular invertible functions on $G$, i.e., the scalar multiples of
characters. This yields an exact sequence $$
0 \rightarrow \chi(G)
\rightarrow k(G)^{(B\times H)}/k^* \rightarrow \bZ \hcD \rightarrow 0,
$$
where $\chi(G)=\chi(C)$ denotes the character group, and $\bZ \hcD$
the free abelian group on $\hcD$. This sequence splits by sending each
$f\in k(G)^{(B\times H)}$ to its restriction to $C$.  Hence, the group
in the middle is canonically $\chi(G)\oplus \bZ\hcD = \chi(C)\oplus
\bZ\hcD$.

On the other hand, assigning to each $f\in k(G)^{(B\times H)}$ its
$H$-weight yields another exact sequence $$
0 \rightarrow N^* =
k(G/H)^{(B)}/k^* \rightarrow \chi(C)\oplus \bZ\hcD = k(G)^{(B\times
  H)}/k^* \rightarrow \chi(H) \rightarrow 0.
$$
Dually, we have a surjective homomorphism $N_C \oplus \bZ \hcD
\onto N$. Its restriction to $N_C$ is injective, since $C$ acts
faithfully. In the analogous sequences for the group $\wH$, the map
$N_C\to \wN$ is zero (since $\wH$ contains $C$), and $\wN_0$ is
generated by the images of colors $\bZ\hcD$.

Denote $\ker(\pi: N\to \wN)$ by $N_{\wH/H}$. Since the action is
assumed to be smart, the map $N_C\to N$ is an isomorphism between
$N_C$ and $N_{\wH/H}$. Changing the action to an equivalent action
composes this isomorphism with an automorphism of $N_C$.

Now, fix one polytope $Q$ in one of the finitely many equivalence
classes modulo $\Aut(N,\pi)$ that we obtained above. This also fixes a
rational polytope $Q'=Q\cap N_{\wH/H}\otimes \bQ$.  The images of
colors belong to $Q$, hence there are only finitely many choices for
the map $\hcD\to Q$. Pick and fix a basis $e_1\dots e_s$ of $N_C$ and
integral points $f_1\dots f_s$ in a fixed multiple $mQ'$ giving a
basis of $N_{\wH/H}$. By switching to an equivalent action, we can
assume that $N_C\to N$ sends $e_i$ to $f_i$. Importantly, if we
replace the polytope $Q$ by another polytope $gQ$, $g\in \Aut(N,\pi)$
in the same equivalence class, then the kernel of $N_C\oplus \bZ\hcD
\onto N$ will not change. Hence, for each of the finitely many
equivalence classes of polytopes we obtained one surjection $\ker
N_C\oplus \bZ\hcD \onto N$ and a single polytope $Q$ in $N$.  Now,
consider the diagram $$
\begin{CD}
  @. 0 \\
  @. @VVV \\
  0 @>>> \wN^* @>>> k(G)^{(B\times H)}/k^*
  @>>> \chi(\wH) @>>> 0 \\
  @. @VVV @| @VVV \\
  0 @>>> N^* @>>> k(G)^{(B\times \wH)}/k^*
  @>>> \chi(H)   @>>> 0 \\
  @. @. @. @VVV \\
  @. @. @. 0 \\
\end{CD}
$$
By what we just have shown, for the cokernel $\chi(H)$ we have only
finitely many choices.  For each of them, the diagram tells us which
characters of $\wH$ vanish on $H\subseteq\wH$, and this determines $H$
since $\wH/H$ is diagonalizable. For each $N^*$, we have finitely many
polytopes in $N$, so finally only finitely many Fano varieties.  The
proof of Theorem~\ref{thm:refined} is now complete.

\end{proof}

\section{Fano invariant of a spherical subgroup} \label{sec:Fano invariant of a spherical subgroup}

By \cite{Brion_Mori}, every spherical $G/H$ has an equivariant Fano
compactification, which is automatically a spherical variety. It is
obtained by running the ``anti Minimal Model Program'', as follows.
Choose an arbitrary projective compactification $X$ of $G/H$. Then by
\cite[Cor.4.7]{Brion_Mori}, there exists a finite sequence $X-\to X'$
of $G$-equivariant divisorial contractions in rays $[C]$ with
$K_X\cdot C>0$ and antiflips such that $-K_{X'}$ is nef; and that
there exists a contraction $X'\to X''$ of a face of $NE(X')$ such that
$X''$ is Fano.

\begin{definition}
  The \emph{Fano invariant} of a spherical subgroup $H\subset G$ is
  $$
  F_G(H) = \max \, \{ \mld(X) \,|\, X \text{ is an equivariant Fano
    compactification of } G/H \}
  $$
\end{definition}
As an immediate consequence of Theorem~\ref{thm:refined}, we have
\begin{corollary}\label{cor:Fano_inv} $F_G(H)$ is well-defined. For
  any $\epsilon>0$, there exist only finitely many $H$ with
  $F_G(H)>\epsilon$, up to an equivalence defined by automorphisms of
  a smart action of a bigger group $G'$. 
\end{corollary} 
By definition, $0<F_G(H) \le 2$, and if there exists a smooth Fano
compactification then $F_G(H)=2$; that is the best case.

\begin{example}
  If $G$ is a torus then for any closed subgroup $H$, $F_G(H)=2$.
  Indeed, in this case $G/H$ is a torus, and can be equivariantly
  compactified to a smooth Fano variety, for example a projective
  space. \end{example}

On the other hand, we have
\begin{lemma}
  If $G$ is a non-commutative reductive group then $\inf_H F_G(H) =0$.
\end{lemma} \begin{proof} We can assume the action is smart.  Let $U$
  be the unipotent radical of a Borel subgroup $B$, and and $S$ be a
  finite subgroup of $T=B\cap \Gss$. By taking $S=T_n$, the subgroup
  of $n$-torsions, we obtain infinitely many spherical subgroups
  $H=US$, and none of them are equivalent. But by
  Corollary~\ref{cor:Fano_inv}, for any $\epsilon>0$ there are only
  finitely many equivalence classes of $H$ with $F_G(H)>~\epsilon$.
\end{proof}

\begin{example}
  If $G$ is a semisimple adjoint group then $F_{G\times G}(\diag
  G)=2$.  Indeed, in this case one has the wonderful compactification
  of De Concini and Procesi \cite{dCP}, which is a smooth Fano.
\end{example}

The proof of the above Lemma may give an impression that perhaps Fano
invariant merely counts the number of connected components of $H$ for
a faithful action of $G$.  But this is not the case. Indeed, for any
spherical subgroup that equals its normalizer $H=\wH$, there exists a
wonderful compactification: in this case the valuation cone $\cV$ is
strictly convex, and one simply takes the toroidal compactification
for the colored fan consisting of $\cV$ and its faces. The wonderful
compactifications are always smooth, and they ``tend to be'' Fanos.
However, $H$ may have many connected components.

\begin{acknowledgements-no}
  The first author was partially supported by NSF grant DMS 0101280.
\end{acknowledgements-no}


\providecommand{\bysame}{\leavevmode\hbox to3em{\hrulefill}\thinspace} \providecommand{\MR}{\relax\ifhmode\unskip\space\fi MR } 
\providecommand{\href}[2]{#2}

\end{document}